\documentclass[12pt,a4paper,reqno]{amsart}
\usepackage{amsmath,amssymb}
\usepackage[left=2.7cm,right=2.7cm,top=2.5cm,bottom=2.5cm]{geometry}

\usepackage{needspace}
\usepackage{xcolor}
\usepackage{hyperref}
\hypersetup{
  colorlinks=false,
  pdfborder={0 0 0.7},
  linkbordercolor={1 0 0},
  citebordercolor={0 1 0}
}

\makeatletter
\newcommand{\partref}[2]{\ref{#1}(#2)}

\let\originalhyperlinkurl\hyper@linkurl
\def\hyper@linkurl#1#2{%
  \begingroup
  \def\@pdfborder{0 0 0}%
  \originalhyperlinkurl{#1}{#2}%
  \endgroup}
\makeatother

\newcommand{\N}{\mathbb{N}}
\newcommand{\Z}{\mathbb{Z}}
\newcommand{\Q}{\mathbb{Q}}
\newcommand{\R}{\mathbb{R}}

\newcommand{\K}{\mathbb{K}}

\DeclareMathOperator{\HF}{HF}
\DeclareMathOperator{\HM}{HM}
\DeclareMathOperator{\SH}{SH}
\DeclareMathOperator{\Gr}{Gr}
\DeclareMathOperator{\Sp}{Sp}
\DeclareMathOperator{\ind}{ind}

\DeclareMathOperator{\supp}{supp}
\DeclareMathOperator{\sign}{sign}

\theoremstyle{plain}
\newtheorem{thm}{Theorem}[section]
\newtheorem{pro}[thm]{Proposition}
\newtheorem{claim}{Claim}
\newtheorem{step}{Step}
\theoremstyle{definition}
\newtheorem{rem}[thm]{Remark}
\numberwithin{equation}{section}
\title[Local homology and closed Reeb orbits in $T^*S^n$]{Local homology and closed Reeb orbits on fiberwise star-shaped hypersurfaces in $T^*S^n$}

\author{Huagui Duan}
\address{School of Mathematics, Jilin University, Changchun, Jinlin, 130012, The People's Republic of China}
\email{duanhg@nankai.edu.cn}

\author[Zihao Qi]{Zihao Qi$^*$}
\address{School of Mathematical Sciences, Nankai University,
  Tianjin 300071, The People's Republic of China}
\email{1120230029@mail.nankai.edu.cn}
\thanks{$^*$ Corresponding author}

\author{Dong Xie}
\address{School of Mathematical Sciences, Nankai University,
  Tianjin 300071, The People's Republic of China;
  Centre de recherches math\'{e}matiques (CRM), Universit\'{e} de Montr\'{e}al,
  Montr\'{e}al, QC H3C 3J7, Canada}
\email{dong.xie@nankai.edu.cn}

\date{2026-12-01}
\keywords{Local Floer homology, closed Reeb orbits, dynamically convex,
  fiberwise star-shaped hypersurface}
\subjclass[2020]{37J46, 37J55, 53D40, 58E05}

\begin{document}

\begin{abstract}
We establish an endpoint property for local Floer homology, showing that
whenever it is nonzero at the lower or upper Conley--Zehnder index, it is one-dimensional and concentrated in that degree. The same property also holds for local equivariant symplectic homology.

As an application, we prove that a dynamically convex closed fiberwise star-shaped
hypersurface in $T^*S^n$, enclosing the zero section
carries at least $\lfloor\frac{n+1}{2}\rfloor$ simple closed Reeb orbits. In addition,
for such a hypersurface with a non-degenerate contact form and finitely many simple closed Reeb orbits, we prove that at least two of these orbits are irrationally elliptic.
\end{abstract}

\maketitle

\section{Introduction and main results}\label{sec:introduction}

\subsection{Background and main results}
Let $M$ be a contact manifold with contact form $\alpha$.
The \textit{Reeb vector field} $R_\alpha$ is defined by
$\alpha(R_\alpha)=1$ and $d\alpha(R_\alpha,\cdot)=0$.
A closed Reeb orbit is a map $\gamma:\R/T\Z\to M$, for some $T>0$,
satisfying $\dot{\gamma}(t)=R_\alpha(\gamma(t))$.
It is called \textit{simple} if it is not a multiple covering of another
closed Reeb orbit.
For $m\in\N$, its $m$-th iteration is
$\gamma^m:\R/mT\Z\to M$, $\gamma^m(t)=\gamma(t)$.
Two simple closed Reeb orbits are \textit{distinct} if they have different
images in $M$.
We omit the word \textit{distinct} when counting simple closed Reeb
orbits.

In this paper, we study closed Reeb orbits on closed fiberwise
star-shaped hypersurfaces in $T^*S^n$, $n\geq3$, enclosing the zero
section, with the contact form induced by the canonical Liouville
$1$-form. The unit cotangent bundle of a Finsler metric on $S^n$ is a
fiberwise star-shaped hypersurface, and its Reeb flow corresponds
to the geodesic flow.

For each $n\geq2$, there are Katok--Ziller metrics on $S^n$
with exactly $2\lfloor\frac{n+1}{2}\rfloor$ prime closed geodesics;
see \cite{Z83}.
It is conjectured that every closed fiberwise star-shaped
hypersurface in $T^*S^n$ enclosing the zero section carries at least
$2\lfloor\frac{n+1}{2}\rfloor$ simple closed Reeb orbits.

For $n=2$, Bangert and Long proved that every Finsler metric
on $S^2$ has at least two distinct prime closed geodesics \cite{BL10}.
Cristofaro-Gardiner, Hryniewicz, Hutchings and Liu subsequently
proved that every Finsler metric on $S^2$ has either two or infinitely
many prime closed geodesics \cite{CHHL26}.
The same alternative holds for fiberwise star-shaped hypersurfaces
in $T^*S^2$: they are diffeomorphic to $\mathbb{R}P^3$, so the
required torsion condition on the first Chern class of the contact
structure is automatic.

For $n\geq3$, the multiplicity conjecture remains open in
general. For Finsler spheres with reversibility $\lambda$ and flag
curvature $K$ satisfying $(\lambda/(1+\lambda))^2<K\leq1$,
Wang proved the sharp lower bound
$2\lfloor\frac{n+1}{2}\rfloor$ for prime closed geodesics under the
bumpy condition \cite{W12}.
Without assuming bumpiness, he obtained at least
$\lfloor\frac{n+1}{2}\rfloor$ prime closed geodesics under the same
curvature condition \cite{W19}.
For $n\geq4$, Duan and Xie obtained at least $n$ prime closed
geodesics under stronger curvature pinching \cite{DX26}.

For contact forms defining the standard contact structure
on the unit cotangent bundle of $S^n$, Abreu and Macarini proved
that dynamical convexity implies the existence of at least two
simple closed Reeb orbits, and that the sharp lower bound
$2\lfloor\frac{n+1}{2}\rfloor$ holds in the non-degenerate case
\cite{AM17}.
They later obtained the same sharp lower bound under weaker index
assumptions, still requiring non-degeneracy \cite{AM26}.
Without assuming non-degeneracy, Ginzburg and G\"urel obtained at
least $\lfloor\frac n2\rfloor-1$ simple closed Reeb orbits on
dynamically convex fiberwise star-shaped hypersurfaces in
$T^*S^n$ \cite{GG20}.

We improve this lower bound using what we call the
\textit{endpoint property of local homology}, which is motivated by \cite[Theorems~1.8 and~4.1]{LLW25}.
For an isolated critical point of a smooth function on a
$d$-dimensional manifold, if a critical group with coefficients
in a field is nonzero in degree $0$ or $d$, then it is
one-dimensional and all critical groups in other degrees vanish
\cite{C93}.
The equivalent statement for local Morse homology is recalled in
\eqref{eq:local-morse-endpoints}.
A similar property holds for the rational $S^1$-equivariant
critical modules of closed geodesics \cite{BL10}.
The following theorem establishes this property for local
Floer homology and local equivariant symplectic homology.

Fix a field $\K$. We use the index conventions
of Section~\ref{sec:index-basic}.

\begin{thm}
\label{endpoint support}
\textup{(i)} Let $(W^{2n},\omega)$ be a symplectic manifold and let
$x$ be an isolated one-periodic orbit of a smooth one-periodic
Hamiltonian $H$, defined on a neighborhood of the graph of $x$
in $S^1\times W$.
Let $\tau$ be any symplectic trivialization of $x^*TW$.
For each $e\in\{\mu_-(x;\tau),\mu_+(x;\tau)\}$,
\begin{equation*}
\dim_{\K}\HF_e(H,x;\K)\leq1.
\end{equation*}
If $\HF_e(H,x;\K)\ne0$, then
$\HF_j(H,x;\K)=0$ for every $j\ne e$.

\textup{(ii)} Let $\alpha$ be a contact form on a cooriented contact
manifold $M^{2n-1}$, $n\geq2$, and let $\gamma=\eta^k$ be an isolated
closed Reeb orbit of $\alpha$, where $\eta$ is simple and
$k\cdot1_{\K}\ne0$.
Let $\xi=\ker\alpha$ and let $\tau$ be any symplectic
trivialization of $\gamma^*\xi$.
For each $e\in\{\mu_-(\gamma;\tau),\mu_+(\gamma;\tau)\}$,
\begin{equation*}
\dim_{\K}\SH^{S^1}_{e}(\gamma;\K)\leq1.
\end{equation*}
If $\SH^{S^1}_e(\gamma;\K)\ne0$, then
$\SH^{S^1}_j(\gamma;\K)=0$ for every $j\ne e$.
\end{thm}
\begin{rem}
\textup{(i)} Both local homology groups are supported in
$[\mu_-(x;\tau),\mu_+(x;\tau)]$ and
$[\mu_-(\gamma;\tau),\mu_+(\gamma;\tau)]$, respectively; see
\eqref{supp Floer} and \eqref{supp}.
Theorem~\ref{endpoint support} states that if either endpoint
belongs to the support, then the homology is one-dimensional
and concentrated in that degree.

\textup{(ii)} Theorem~\ref{endpoint support} removes the
additional assumption on $\mu_+$ in
\cite[Theorems~1.8 and~4.1]{LLW25}.
A related argument appears in the proof of
\cite[Theorem~D]{CGG26}.

\end{rem}

We now apply Theorem~\ref{endpoint support} to improve the
multiplicity bound obtained in \cite[Theorem~1.6]{GG20}.
A contact form $\alpha$ on a fiberwise star-shaped hypersurface in
$T^*S^n$ is \textit{dynamically convex} if $\mu_-(\gamma)\geq n-1$
for every closed Reeb orbit $\gamma$.
For a Finsler metric on $S^n$, the pinching condition
$(\lambda/(1+\lambda))^2<K\leq1$ implies dynamical convexity
of the induced contact form on its unit cotangent bundle
\cite{AM17}.

\begin{thm}
\label{multiplicity}
Let $M^{2n-1}$ be a closed fiberwise star-shaped hypersurface
in $T^*S^n$ enclosing the zero section, where $n\geq3$.
Let $\alpha$ be the restriction of the canonical Liouville
$1$-form to $M$.
If $\alpha$ is dynamically convex, then $M$ carries
at least $\lfloor\frac{n+1}{2}\rfloor$ simple closed Reeb orbits.
\end{thm}

\begin{rem}(i) For $n\geq5$, this improves both general lower bounds
$2$ and $\lfloor\frac n2\rfloor-1$ recalled above.
It also extends the lower bound in \cite{W19} from Finsler spheres
satisfying $(\lambda/(1+\lambda))^2<K\leq1$ to dynamically convex
fiberwise star-shaped hypersurfaces in $T^*S^n$.

(ii) Restricted contact type hypersurfaces are not considered here, since generally its positive equivariant symplectic homology may not be isomorphic to the unit cotangent bundle, see Remark \ref{rem: restricted contact type} for details.
\end{rem}

A closed Reeb orbit is \textit{elliptic} if all eigenvalues of its
linearized Poincar\'{e} map $P$ lie on the unit circle.
It is \textit{irrationally elliptic} if $P$ can be represented,
in a symplectic basis, by a direct sum of $2\times2$ rotation
matrices with rotation angles $\theta_j$ satisfying
$\theta_j/\pi\notin\Q$.
The orbit is \textit{non-degenerate} if $1\notin\sigma(P)$.
A contact form is \textit{non-degenerate} if all its closed Reeb
orbits are non-degenerate.

For fiberwise star-shaped hypersurfaces in $T^*S^2$,
finiteness of the number of simple closed Reeb orbits implies
that there are exactly two and both are irrationally elliptic
\cite{CHHL26}.
For higher-dimensional Finsler spheres, the existence of at least
two elliptic prime closed geodesics was proved under curvature
pinching and the assumption of finitely many prime closed geodesics
in \cite{D16}.
Hamid and Wang extended this result to compact Finsler manifolds.
Under the additional bumpy assumption, they showed that the two
closed geodesics can be chosen to be irrationally elliptic
\cite{HW22}.

For dynamically convex Reeb flows on closed star-shaped
hypersurfaces in $\R^{2n}$, $n\geq2$, the existence of at least two
irrationally elliptic simple closed orbits was proved in
\cite{LLW25} under the assumption that there are
exactly $n$ simple closed orbits, and in
\cite{CGG26} under the assumption that there are
finitely many simple closed orbits.

For fiberwise star-shaped hypersurfaces in $T^*S^n$,
we obtain the following result.

\begin{thm}
\label{irrationally elliptic}
Let $M^{2n-1}$ be a closed fiberwise star-shaped hypersurface
in $T^*S^n$ enclosing the zero section, where $n\geq3$.
Let $\alpha$ be the restriction of the canonical Liouville
$1$-form to $M$.
If $\alpha$ is dynamically convex and non-degenerate and has
only finitely many simple closed Reeb orbits, then at least
two of them are irrationally elliptic.
\end{thm}

A Reeb flow with only finitely many simple closed orbits is called a
\textit{pseudo-rotation}.
The Katok--Ziller examples give dynamically convex pseudo-rotations
whose closed orbits are all irrationally elliptic \cite{Z83,AM17}.
These examples motivate the conjecture that all closed Reeb orbits
of a pseudo-rotation on a closed fiberwise star-shaped hypersurface
in $T^*S^n$ enclosing the zero section are irrationally elliptic.
Theorem~\ref{irrationally elliptic} is a step toward this conjecture
under dynamical convexity and non-degeneracy.

We expect that, without the non-degeneracy assumption, the
remaining hypotheses of Theorem~\ref{irrationally elliptic}
still guarantee at least two simple elliptic closed Reeb orbits.
One obstacle is that it remains open whether a
symplectically degenerate maximum (SDM) that is not simple
forces infinitely many simple closed Reeb orbits in this setting.
Fernandes and Pering proved this for isolated SDMs on the
boundary of a Liouville domain $W^{2n}$, $n\geq2$, with
$c_1(TW)=0$, provided that the SDM represents a non-torsion
class in $H_1(W;\Z)$ \cite[Theorem~2]{FP26}.
This condition fails here, since the domain bounded by $M$
retracts onto $S^n$ and hence has trivial first homology.

\subsection{Outline of the proofs}

For (i) of Theorem~\ref{endpoint support}, we use generating
functions and the splitting lemma to identify local Floer
homology, up to a grading shift, with the local Morse homology
of a function on a $\nu(x)$-dimensional space, where $\nu(x)$
is the nullity \cite{G10,GG10}.
This identification proceeds through the local Lagrangian Floer homology.
We determine this shift and show that $\mu_-(x;\tau)$ and
$\mu_+(x;\tau)$ correspond to degrees $0$ and $\nu(x)$,
respectively.
The endpoint property for local Morse homology then gives
the result.
Part~\textup{(ii)} follows by comparing local equivariant
symplectic homology with the local Floer homology of the
return map.

The proof of Theorem~\ref{multiplicity} combines positive
equivariant symplectic homology with index iteration theory.
We first obtain $\lfloor\frac n2\rfloor-1$ simple closed Reeb
orbits by the index iteration argument in \cite{LZ02,GG20}.
The Morse inequalities and index estimates then yield an
iterate with nonzero local homology at its upper
Conley--Zehnder index.
By the endpoint property, this homology is concentrated in
that degree, which shows that the underlying simple orbit
is distinct from those already found.
When $n$ is odd, an analysis of the equality cases in the
index estimates shows that, if there were only
$\lfloor\frac n2\rfloor$ simple closed Reeb orbits, one of
them would be an SDM.
Proposition~\ref{prop:sdm} would then imply infinitely many
simple closed Reeb orbits, a contradiction.

For Theorem~\ref{irrationally elliptic}, we start with the
orbit obtained in Step~\ref{step:2} of the proof of
Theorem~\ref{multiplicity}.
Under non-degeneracy, the index estimates and the basic
normal form of its linearized Poincar\'{e} map show that
this orbit is irrationally elliptic.
The symmetry in the common index jump theorem then yields
a second irrationally elliptic orbit.

The proof of Theorem~\ref{endpoint support}
is given in Section~\ref{sec:proof-endpoint}.
The proofs of Theorems~\ref{multiplicity} and~\ref{irrationally elliptic}
are given in Sections~\ref{sec:proof-multiplicity}
and~\ref{sec:proof-ellipticity}, respectively.

\medskip

\textbf{Some explanations and AI disclosure.} This paper is actually the revised version of the manuscript with arXiv:2603.07894 in March 2026, according to feedbacks and suggestions from several experts. In this version, ChatGPT-6-Astra was only used for polishing the introduction and searching references in the proof of Theorem 1.1.

\section{Local homology and the endpoint property}\label{sec:floer}

\subsection{Conley--Zehnder indices}\label{sec:index-basic}

For $d\in\N$ and $T>0$, let
\[
  \mathcal{P}_T(2d)
  =\{\Phi\in C([0,T],\Sp(2d))\mid \Phi(0)=I_{2d}\}.
\]
The nullity of $\Phi\in\mathcal{P}_T(2d)$ is
$\nu(\Phi)=\dim\ker(\Phi(T)-I_{2d})$.
The path is non-degenerate if $\nu(\Phi)=0$.
For a non-degenerate path, we denote its Conley--Zehnder
index by $\mu_{\mathrm{CZ}}(\Phi)$, normalized by
$\mu_{\mathrm{CZ}}(t\mapsto e^{2\pi i\lambda t})=1$
for $t\in[0,1]$ and $0<\lambda<1$.

The lower and upper semicontinuous extensions of this index
to $\mathcal{P}_T(2d)$ are
\[
  \mu_-(\Phi)=\liminf_{\Psi\to\Phi}\mu_{\mathrm{CZ}}(\Psi),
  \qquad
  \mu_+(\Phi)=\limsup_{\Psi\to\Phi}\mu_{\mathrm{CZ}}(\Psi),
\]
where $\Psi\in\mathcal{P}_T(2d)$ ranges over non-degenerate
paths converging uniformly to $\Phi$.
The difference between them are
\begin{equation}
  \mu_+(\Phi)-\mu_-(\Phi)=\nu(\Phi),
  \label{mu+-mu-}
\end{equation}
see \cite[Section~6.1]{L02}.
For non-degenerate $\Phi$, both indices equal
$\mu_{\mathrm{CZ}}(\Phi)$.
We also use the Robbin--Salamon index
$\mu_{\mathrm{RS}}(\Phi)$, which agrees with
$\mu_{\mathrm{CZ}}(\Phi)$ on non-degenerate paths
\cite{RS93}.

For $\Phi\in\mathcal{P}_T(2d)$ and $m\in\N$, the
$m$-th iteration $\Phi^m\in\mathcal{P}_{mT}(2d)$ is defined by
\begin{equation}
  \Phi^m(t)=\Phi(t-jT)\Phi(T)^j,
  \qquad
  jT\le t\le(j+1)T,\quad j=0,\ldots,m-1.
  \label{3.9}
\end{equation}
The mean index is
$\hat{\mu}(\Phi)=\lim_{m\to\infty}\mu_-(\Phi^m)/m$.
It satisfies $\hat{\mu}(\Phi^m)=m\hat{\mu}(\Phi)$.
Moreover, it is related to the upper and lower extensions by
\begin{equation}
  |\mu_\pm(\Phi)-\hat{\mu}(\Phi)|\le d.
  \label{mu-hat{mu}}
\end{equation}

We use the convention $\omega(X_{H_t},\cdot)=-dH_t$
for Hamiltonian vector fields.
For a closed Hamiltonian orbit $x$ on $(W,\omega)$,
choose a symplectic trivialization $\tau$ of $x^*TW$.
For a closed Reeb orbit $\gamma$ on $(M,\alpha)$, choose
$\tau$ on $\gamma^*\xi$, where $\xi=\ker\alpha$.
The linearized flow, restricted to $\xi$ in the Reeb case,
then determines a symplectic path.
The indices and nullity are denoted by
$\mu_{\mathrm{RS}}(x;\tau)$, $\mu_\pm(x;\tau)$,
$\hat{\mu}(x;\tau)$ and $\nu(x)$ in the Hamiltonian case,
and similarly with $x$ replaced by $\gamma$ in the Reeb case.

The indices depend only on the homotopy class of $\tau$,
whereas the nullity is independent of $\tau$.
We use the same trivialization to grade the local homology
groups.
Changing $\tau$ shifts the indices by an even integer.
For iterated orbits $x^m$ and $\gamma^m$, we use the iterated
trivializations $\tau^m$.

For a Hamiltonian orbit contractible in $W$, or a Reeb orbit
contractible in a Liouville filling $W$, we use a
trivialization induced by a capping disk.
If $c_1(TW)|_{\pi_2(W)}=0$, its homotopy class is independent
of the capping disk, and we omit $\tau$ from the notation.
The construction for Reeb orbits is recalled in
Section~\ref{sec:global-equivariant-homology}.

\Needspace{6\baselineskip}

\subsection{Local homology}

\subsubsection{Local Morse homology and local Floer homology}
\label{sec:local-homology}

Let $f$ be a smooth function on a manifold $M^m$, and let $p$
be an isolated critical point.
Choose a sufficiently small relatively compact neighborhood
$U$ of $p$ such that $\overline U$ contains no other critical
points of $f$.
For a sufficiently $C^1$-small Morse perturbation $\tilde f$
and a Morse--Smale metric, the critical points of $\tilde f$
in $U$ generate a local Morse complex over a field $\K$.
Its differential counts negative gradient trajectories
contained in $U$, with signs determined by orientations
of the unstable manifolds.
The resulting local Morse homology $\HM_*(f,p)$
is independent of these choices;
see \cite{G10,GG10}.

The local Morse homology $\HM_*(f,p)$ is isomorphic,
with the same grading, to the critical groups of $f$ at $p$
with coefficients in $\K$.
Both are computed by the relative homology of a Gromoll--Meyer
pair; see \cite[Theorem~7.1]{RV14} and
\cite[Chapter~I, Theorem~5.2]{C93}.

For $e\in\{0,m\}$, if $\HM_e(f,p)\ne0$, then
\begin{equation}\label{eq:local-morse-endpoints}
  \HM_j(f,p)\cong
  \begin{cases}
    \K, & j=e,\\
    0, & j\ne e,
  \end{cases}
  \qquad j\in\mathbb{Z};
\end{equation}
see \cite[Chapter~I, Corollary~5.1]{C93}.

Let $(W^{2n},\omega)$ be a symplectic manifold, and let $x$
be an isolated one-periodic orbit of a smooth one-periodic
Hamiltonian $H$ defined near the graph of $x$ in $S^1\times W$.
Write $\varphi_t=\varphi_H^t$ and set $p=x(0)$.
Let $\tau$ be any symplectic trivialization of $x^*TW$.
The linearized Hamiltonian flow is represented by
\[
  \Phi_\tau(t)
  =\tau_t\circ d\varphi_t(p)\circ\tau_0^{-1},
  \qquad t\in[0,1],
\]
where $\tau_t:T_{x(t)}W\to\mathbb{R}^{2n}$.

Choose a sufficiently small isolating tubular neighborhood
$V$ of the graph of $x$ in $S^1\times W$, and extend $\tau$
to a symplectic trivialization of the pullback of $TW$
under the projection $V\to W$.
Choose a smooth one-periodic perturbation $\tilde H$ of $H$,
sufficiently $C^2$-close to $H$ and equal to $H$ outside $V$,
whose one-periodic orbits with graphs in $V$ are non-degenerate.
These orbits generate the local Floer complex over $\K$,
graded by their Conley--Zehnder indices.
For a regular one-periodic family $J_t$ of
$\omega$-compatible almost complex structures,
the differential counts Floer cylinders
$u:\mathbb{R}\times S^1\to W$ satisfying
$(t,u(s,t))\in V$ for all $(s,t)$,
with signs determined by coherent orientations.
The resulting local Floer homology $\HF_*(H,x)$ is
independent of $V$, $\tilde H$, and $J_t$;
see \cite{G10,GG10}.

For $\tilde H$ sufficiently $C^2$-close to $H$, the
Conley--Zehnder indices of the generators lie in
$[\mu_-(x;\tau),\mu_+(x;\tau)]$.
Hence
\begin{equation}\label{supp Floer}
  \supp\HF(H,x)
  \subset[\mu_-(x;\tau),\mu_+(x;\tau)].
\end{equation}

\subsubsection{Local equivariant symplectic homology}
\label{sec:equivariant-homology}

Let \(\gamma=\eta^k\) be an isolated closed Reeb orbit on
$(M^{2n-1},\alpha)$ with period \(T\),
where $\eta$ is simple.
We parametrize $\gamma$ on $S^1=\mathbb{R}/\mathbb{Z}$, so that
$\dot \gamma=T R_\alpha(\gamma)$.
Choose a tubular neighborhood $U$ of $\gamma(S^1)$ containing
no other closed Reeb orbit with period sufficiently close
to $T$, up to time shift.
On
\(\bigl((1-\varepsilon,1+\varepsilon)\times U,d(r\alpha)\bigr)\),
take \(H=h(r)\) with \(h'(1)=T\) and \(h''(1)>0\).
Since $X_H=h'(r)R_\alpha$, for sufficiently small
$\varepsilon$ the Hamiltonian $H$ has an isolated circle
of one-periodic orbits
\[
  S_\gamma=\{t\mapsto(1,\gamma(t+\theta)):\theta\in S^1\}.
\]
Let $\tau$ be any symplectic trivialization of $\gamma^*\xi$,
where $\xi=\ker\alpha$, and extend it along $(1,\gamma)$
by the frame $(\partial_r,R_\alpha)$.

Define the local equivariant symplectic homology by
\[
  \SH^{S^1}_*(\gamma;\K)
  :=\HF^{S^1}_*(H,S_\gamma;\K),
\]
where the right-hand side is the local
\(S^1\)-equivariant Floer homology of \(H\) at \(S_\gamma\)
\cite[Section~2.3.1]{GG20}.
For fixed $\tau$, these graded groups are independent of
the auxiliary choices.

Set $p=\gamma(0)$ and choose a transverse disk
$(B^{2n-2},d\alpha|_B)$ through $p$ with $T_pB=\xi_p$,
meeting $\eta$ only at $p$.
Let $\varphi_\eta$ be the first-return map on $B$ near $p$,
and set $\varphi_\gamma=\varphi_\eta^k$.
The isolation of $\gamma$ implies that $p$ is an isolated
fixed point of $\varphi_\gamma$.
Choose a local one-periodic Hamiltonian $H_\gamma$ on $B$
generating $\varphi_\gamma$ near $p$ as in \cite[Lemma~5.2]{HM15}.
The Hamiltonian flow fixes $p$ and linearized Hamiltonian path under the trivialization $\tau_\gamma$ of $S^1\times T_pB$ induced by $\tau$ is a reparametrization of the linearized Reeb path of $\gamma$ in $\tau$.
Thus $\mu_\pm(p;\tau_\gamma)=\mu_\pm(\gamma;\tau)$.

The projections of periodic orbits close to $S_\gamma$
are freely homotopic to $\gamma$ in $U$.
Let $\pi:U_k\to U$ be the connected $k$-fold cyclic cover,
equipped with $\pi^*\alpha$.
The orbit $\gamma$ lifts to a simple isolated Reeb orbit
$\widetilde \gamma$, whose first-return map on the lifted
transverse disk is $\varphi_\gamma=\varphi_\eta^k$.
Equip $\widetilde \gamma$ with the trivialization induced by $\tau$.
Assume that $k\cdot1_{\K}\ne0$.
By \cite[Lemma~2.9]{CGG25} and
\cite[Theorem~1.2]{F20}, respectively,
\[
  \begin{aligned}
    \SH^{S^1}_j(\gamma;\K)
    &\cong\SH^{S^1}_j(\widetilde \gamma;\K)^{\mathbb{Z}_{k}},\\
    \SH^{S^1}_j(\widetilde \gamma;\K)
    &\cong\HF_j(H_\gamma,p;\K).
  \end{aligned}
\]
The superscript $\mathbb{Z}_{k}$ denotes invariants under
the action induced by deck transformations.

Taking dimensions gives
\begin{equation}\label{eq:equivariant-sh-comparison}
  \dim_{\K}\SH^{S^1}_j(\gamma;\K)
  \leq \dim_{\K}\HF_j(H_\gamma,p;\K)
  \qquad\text{for every }j\in\Z.
\end{equation}
For $k=1$, the two groups are isomorphic in each degree.
Together with \eqref{supp Floer}, this gives
\begin{equation}\label{supp}
  \supp\SH^{S^1}(\gamma;\K)\subset[\mu_-(\gamma;\tau),\mu_+(\gamma;\tau)]
  \qquad\text{if }k\cdot1_{\K}\ne0.
\end{equation}

A non-degenerate orbit $\gamma=\eta^k$ is called \textit{good}
if its Conley--Zehnder index has the same parity as that of
$\eta$, and \textit{bad} otherwise.
Then $\SH^{S^1}_*(\gamma;\mathbb{Q})$ is one-dimensional
and concentrated in degree $\mu_-(\gamma;\tau)=\mu_+(\gamma;\tau)$
if $\gamma$ is good, and vanishes if $\gamma$ is bad;
see \cite[Example~3.6]{GG20}.

\subsection{Proof of Theorem~\ref{endpoint support}}
\label{sec:proof-endpoint}

(i) Set $p=x(0)$, $\varphi_t=\varphi_H^t$, and
$\varphi=\varphi_H^1$.
Let $\Phi_\tau$ be the symplectic path defined in
Section~\ref{sec:local-homology}.
Choose a Darboux chart centered at $p$ whose differential
at $p$ is $\tau_0$.
In this chart, $\varphi$ is a local Hamiltonian
diffeomorphism of $\mathbb{R}^{2n}$ fixing $0$.

Let $\Delta$ be the diagonal of
$(\mathbb{R}^{2n},-\omega_{std})\times
(\mathbb{R}^{2n},\omega_{std})$, and choose a common
Lagrangian complement $N$ to $\Delta$ and
$\Gr(d\varphi_p)$.
Then this symplectic vector space can be identified with
\[
  \left(T^*\Delta,\sum_{i=1}^{2n}dq_i\wedge d\zeta_i\right),
\]
where $q$ is the base coordinate and $\zeta$ is the fiber
coordinate.
Over a sufficiently small ball $B\subset\Delta$
centered at $0$, the diagonal is the zero section
$0_B$ of $T^*B$, and
$\Gr(\varphi)=\Gr(dF)$ for a smooth function
$F:B\to\mathbb{R}$ with $F(0)=0$,
called the generating function of $\varphi$.
We refer to \cite[Section~4.2]{GG10} for generating
functions when the Hamiltonian diffeomorphism is not
$C^1$-close to the identity.
Since $p$ is an isolated fixed point of $\varphi$,
the origin is an isolated critical point of $F$.

We use local Lagrangian Floer homology for the pair
$(\Delta,\Gr(\varphi))$ at $(p,p)$.
In a sufficiently small neighborhood of $(p,p)$, the two
Lagrangian branches are exact disks.
Fix spin structures on these disks and use the grading determined
by $t\mapsto\Gr(\Phi_\tau(t))$ relative to the diagonal.
After a sufficiently small Hamiltonian perturbation making the
intersections transverse, the intersection points near $(p,p)$
generate the local Floer complex.
For a regular almost complex structure, its differential counts
rigid Floer strips contained in the neighborhood, with boundary
on the two Lagrangian branches.
The signs are determined by coherent orientations;
see \cite[Section~4.3]{M00}.
We denote the homology by
$\HF_*(\Delta,\Gr(\varphi),(p,p);\K)$.
The homology is independent of the sufficiently small perturbation
and the regular almost complex structure.
For the local construction and the continuation argument,
see \cite[Section~3]{S99}.

Write a Floer cylinder $u$ in the coordinates
$v(s,t)=\varphi_t^{-1}(u(s,t))$.
The periodicity of $u$ becomes
$v(s,0)=\varphi(v(s,1))$.
Folding at $t=\tfrac12$ gives strips with boundary
on $\Delta$ and $\Gr(\varphi)$.
Linearizing this change of variables identifies the
determinant lines of the corresponding Floer operators.
Choose coherent orientations under these identifications
so that corresponding rigid trajectories have the same sign.
With the induced Floer data, the change of variables followed
by folding identifies the local Hamiltonian and Lagrangian Floer complexes.
Hence
\[
  \HF_*(H,x;\K)
  \cong
  \HF_*(\Delta,\Gr(\varphi),(p,p);\K).
\]
For each non-degenerate orbit of a small perturbation,
let $\widetilde\Phi$ be its linearized symplectic path
in the induced trivialization.
The Maslov index of $t\mapsto\Gr(\widetilde\Phi(t))$
relative to $\Delta$ equals
$\mu_{\mathrm{CZ}}(\widetilde\Phi)$;
see \cite[Remark~5.4]{RS93}.
Thus this isomorphism preserves the grading.

In the coordinates $(q,\zeta)$,
$T_{(p,p)}\Delta$ and $T_{(p,p)}\Gr(\varphi)$
correspond to $T_0B\oplus\{0\}$ and
$\Gr(d^2F_0)$, respectively.
Since
\[
  (T_0B\oplus\{0\})\cap\Gr(d^2F_0)
  =\ker(d^2F_0)\oplus\{0\},
\]
we obtain
\[
  \nu:=\dim\ker d^2F_0=\dim\ker(d\varphi_p-I).
\]
By \cite[Section~4.1.2]{GG20} and \eqref{mu+-mu-},
\begin{equation}\label{eq:graph-nullity}
  \mu_\pm(x;\tau)
  =\mu_{\mathrm{RS}}(\Phi_\tau)\pm\frac{\nu}{2}.
\end{equation}

Choose an isolating neighborhood $U\subset\subset B$
of the critical point $0$ of $F$, and a Morse function
$G$ sufficiently $C^2$-close to $F$.
Set $F_\lambda=(1-\lambda)F+\lambda G$ for
$0\leq\lambda\leq1$.
For a critical point $z\in U$ of $G$, the tangent spaces
to $\Gr(dF_\lambda)$ at
$(\lambda z,dF_\lambda(\lambda z))$ form the path
\[
  \lambda\longmapsto\Gr(d^2(F_\lambda)_{\lambda z}),
\]
which joins $\Gr(d^2F_0)$ to $\Gr(d^2G_z)$.
To compute the Floer degree $|z|$, append this path
to $t\mapsto\Gr(\Phi_\tau(t))$, expressed in
the coordinates $(q,\zeta)$.
Additivity and the localization formula
\cite[Theorem~2.3 and Remark~5.4]{RS93} give
\begin{equation}\label{eq:terminal-degree}
  \begin{aligned}
    |z|
    &=\mu_{\mathrm{RS}}(\Phi_\tau)
      +\frac12\bigl(\sign(d^2G_z)-\sign(d^2F_0)\bigr)\\
    &=\mu_-(x;\tau)
      +\ind_{-G}(z)-\ind_{-F}(0).
  \end{aligned}
\end{equation}
Here $\ind_h(z)$ denotes the number of negative eigenvalues
of $d^2h_z$, counted with multiplicity.
The second equality in \eqref{eq:terminal-degree}
follows from \eqref{eq:graph-nullity} and the identities
\[
  \sign(d^2G_z)=2\ind_{-G}(z)-2n,
  \qquad
  \sign(d^2F_0)=2\ind_{-F}(0)-(2n-\nu).
\]

Define $T_G:T^*B\to T^*B$ by
\[
  T_G(q,\zeta)=(q,\zeta-dG_q).
\]
Then $T_G$ is symplectic and satisfies
$T_G(0_B)=\Gr(-dG)$ and
$T_G(\Gr(dG))=0_B$.
Choose a Riemannian metric $g$ such that $(-G,g)$
is Morse--Smale, and rescale $g$ to satisfy the Hessian
bound in \cite[(39)]{M00}.
Use the almost complex structures from
\cite[proof of Proposition~27]{M00} for
$(\Gr(-dG),0_B)$, pulled back by $T_G$.
With this choice, negative gradient trajectories
$\ell$ of $-G$ correspond bijectively to Floer strips
$w_\ell$ for $(0_B,\Gr(dG))$, with
\begin{equation}\label{eq:floer-morse-trajectories}
  w_\ell(s,t)
  =(\ell(s),t\,dG_{\ell(s)}),
  \qquad
  \dot\ell(s)=-\nabla_g(-G)(\ell(s)).
\end{equation}
By \eqref{eq:floer-morse-trajectories},
$w_\ell(\mathbb{R}\times[0,1])\subset T^*U$
if and only if $\ell(\mathbb{R})\subset U$.
Choose coherent orientations as in
\cite[Proposition~22 and Remark~24]{M00} so that
the correspondence in \eqref{eq:floer-morse-trajectories}
preserves orientations.
The local Floer differential for $(0_B,\Gr(dG))$
then agrees with the local Morse differential of $-G$
over $\K$.

Choose a complement $E$ to $\ker d^2F_0$ in $T_0B$.
By the splitting lemma
\cite[Chapter~I, Theorem~5.1]{C93},
there are local coordinates in which
\[
  F(u,v)=Q_E(u)+f(v),
  \qquad
  u\in E,\quad v\in\ker d^2F_0,
\]
where $Q_E$ is a non-degenerate quadratic form on $E$,
$\ind_{-Q_E}(0)=\ind_{-F}(0)$, and $d^2f_0=0$.
The origin is an isolated critical point of $f$.
Therefore,
\begin{equation}\label{eq:graph-Morse-reduced}
  \HF_j(H,x;\K)
  \cong
  \HM_{j-\mu_-(x;\tau)+\ind_{-F}(0)}(-F,0;\K)
  \cong
  \HM_{j-\mu_-(x;\tau)}(-f,0;\K).
\end{equation}
The local Morse complex of $-G$ in $U$ computes
$\HM_*(-F,0;\K)$.
The correspondence in \eqref{eq:floer-morse-trajectories},
together with the grading formula
\eqref{eq:terminal-degree}, therefore gives the first
isomorphism in \eqref{eq:graph-Morse-reduced}.
The second follows from the identification with critical
groups in Section~\ref{sec:local-homology} and the shifting theorem
applied to $-F=-Q_E-f$; see
\cite[Chapter~I, Theorem~5.4]{C93}.

Since $\dim\ker d^2F_0=\nu=\mu_+(x;\tau)-\mu_-(x;\tau)$,
part (i) follows from \eqref{eq:graph-Morse-reduced}
and the endpoint property \eqref{eq:local-morse-endpoints}.

\Needspace{5\baselineskip}
(ii) Choose $H_\gamma$ and $\tau_\gamma$ as in
Section~\ref{sec:equivariant-homology}.
Let
$e\in\{\mu_-(\gamma;\tau),\mu_+(\gamma;\tau)\}$.
By \eqref{eq:equivariant-sh-comparison} and part (i),
\[
  \dim_{\K}\SH^{S^1}_e(\gamma;\K)
  \leq
  \dim_{\K}\HF_e(H_\gamma,z_\gamma;\K)
  \leq1.
\]
If $\SH^{S^1}_e(\gamma;\K)\ne0$, then
$\HF_e(H_\gamma,z_\gamma;\K)\ne0$.
Part (i) and \eqref{eq:equivariant-sh-comparison} then give
$\SH^{S^1}_j(\gamma;\K)=0$ for every $j\ne e$.
This proves (ii).

\Needspace{9\baselineskip}
\section{Positive equivariant symplectic homology}
\label{sec:global-equivariant-homology}

\subsection{Definition and computation}
Let $W\subset T^*S^n$ be
a compact fiberwise star-shaped domain containing the zero section,
where $n\geq3$.
Let $\lambda$ be the canonical Liouville form, and set
$M=\partial W$ and $\alpha=\lambda|_M$.
Then $(W,d\lambda)$ is a Liouville domain with contact boundary
$(M,\alpha)$.
Every closed Reeb orbit of $\alpha$ is contractible in $W$.
Moreover,
\begin{equation}\label{eq:c1-vanishing}
  c_1(TW)|_{\pi_2(W)}=0.
\end{equation}

Let $Z$ be the Liouville vector field, defined by
$\iota_Zd\lambda=\lambda$.
Along $M$, the symplectic splitting is
\[
  TW|_M=\xi\oplus\operatorname{span}\{Z,R_\alpha\},
  \qquad \xi=\ker\alpha.
\]
For a closed Reeb orbit $\gamma$, a capping trivialization of
$\gamma^*\xi$, together with the frame $(Z,R_\alpha)$, extends over
a capping disk in $W$.

We write $\SH^{S^1,+}_*(W)$ for the positive
$S^1$-equivariant symplectic homology of $W$ in the contractible
free homotopy class, with coefficients in $\mathbb{Q}$.
We grade this homology by the Conley--Zehnder index.
It is obtained as a direct limit of filtered $S^1$-equivariant
Floer homology groups for admissible Hamiltonians on the symplectic
completion of $W$; see \cite[Section~3.1.1]{GG20}.
The superscript $+$ indicates that the contribution of constant
Hamiltonian orbits is excluded.

By \cite[Proposition~3.12]{GG20},
\begin{equation}
  \dim_{\mathbb{Q}}\SH^{S^1,+}_k(W) = \left\{\begin{array}{ll}
    0 \quad if\quad k=n\ (mod\ 2)\ or\ k<n-1,  \\
    2 \quad if\quad k=j(n-1)\ for\ j>1\ when\ n\ is\ odd, \\
    2 \quad if\quad k=j(n-1)\ for\ odd\ j>1\ when\ n\ is\ even, \\
    1 \quad otherwise. \end{array}\right. \label{b_k}
\end{equation}

For an isolated closed Reeb orbit $\gamma$, we write
$\SH^{S^1}_*(\gamma)=\SH^{S^1}_*(\gamma;\mathbb{Q})$ for the local groups
from Section~\ref{sec:equivariant-homology}.
Let $\mathcal P(\alpha)$ denote the set of closed Reeb
orbits of $\alpha$.
For $\gamma:\R/T\Z\to M$ in $\mathcal P(\alpha)$, its action is
\[
  \mathcal A(\gamma):=\int_\gamma\alpha=T.
\]

\Needspace{7\baselineskip}
\begin{pro}
\label{prop:cotangent-homology}
\textnormal{(\cite[Corollary~3.14]{GG20})}
Assume that all closed Reeb orbits of $\alpha$ are
isolated. For every $j\in\mathbb N$, there exist closed Reeb orbits
$\eta_0,\ldots,\eta_{n-1}$ such that
\begin{equation}\label{y_i}
  \SH^{S^1}_{2i+(2j-1)(n-1)}(\eta_i)\neq0,
  \qquad i=0,\ldots,n-1,
\end{equation}
and their actions satisfy
\begin{equation}
\mathcal{A}(\eta_0)<\mathcal{A}(\eta_1)<\cdots<\mathcal{A}(\eta_{n-1}).\label{action<}
\end{equation}
\end{pro}
\begin{rem}\label{rem: restricted contact type}
For general restricted contact type hypersurfaces, the transfer map between $\lambda W$ and $W$ may not be isomorphic, so the argument in \cite[Proposition~3.12]{GG20} does not follow and $b_k$ can be different from (\ref{b_k}). However, Proposition \ref{prop:cotangent-homology} still holds.
\end{rem}

\subsection{Resonance identity and Morse inequalities}

\leavevmode Let $\gamma\in\mathcal P(\alpha)$ be a closed Reeb orbit
all of whose iterates are isolated.
Its local Euler characteristic is
\begin{equation}
\chi(\gamma)= \sum_{m\in\Z}(-1)^m\dim{\SH^{S^1}_m(\gamma)},
\end{equation}
where the sum is finite by \eqref{supp}.
When $\gamma$ is non-degenerate,
$\chi(\gamma)=(-1)^{\mu_-(\gamma)}$ if $\gamma$ is good,
and $\chi(\gamma)=0$ if it is bad.
The local mean Euler characteristic of $\gamma$ is
\[
  \hat\chi(\gamma)
  =\lim_{j\to\infty}\frac{1}{j}\sum_{k=1}^j\chi(\gamma^k).
\]
By \cite[Theorem~2.8 and Remark~3.13]{GG15}, each sequence
$\{\chi(\gamma_i^k)\}_{k\geq1}$ is periodic.
For each $i$, choose $L_i\in\N$ such that
$\chi(\gamma_i^{k+L_i})=\chi(\gamma_i^k)$ for every $k\geq1$.
Then
\begin{equation}\label{hat chi finite sum}
  \hat\chi(\gamma_i)
  =\frac{1}{L_i}\sum_{k=1}^{L_i}\chi(\gamma_i^k)
\end{equation}is finite.

Assume now that $\alpha$ has finitely many simple closed Reeb
orbits $\gamma_1,\ldots,\gamma_r$, with $\hat\mu(\gamma_i)>0$ for every $i$.
Then every iterate $\gamma_i^k$ is isolated.
Set $b_k=\dim_{\mathbb{Q}}\SH^{S^1,+}_k(W)$.
The positive mean Euler characteristic is
\begin{equation}
      \chi_+
  :=\lim_{j\to\infty}\frac{1}{j}\sum_{k=1}^j(-1)^k b_k
  =\begin{cases}
    -\dfrac{n}{2n-2}, & n\text{ is even},\\[3pt]
    \dfrac{n+1}{2n-2}, & n\text{ is odd},
  \end{cases}\label{chi_+}
\end{equation}

where the value follows from \eqref{b_k}.
By the resonance identity (see \cite[Theorem~1.5]{HM15}
and \cite[Section~3.5]{GGM18}),
\begin{equation}
  \sum_{i=1}^r\frac{\hat{\chi}(\gamma_i)}{\hat{\mu}(\gamma_i)}=\chi_+.
  \label{MII}
\end{equation}

Set
\begin{equation*}
  c_m=\sum_{i=1}^r\sum_{k=1}^\infty\dim{\SH^{S^1}_m(\gamma_i^k)}.
\end{equation*}
By \eqref{supp} and \eqref{mu-hat{mu}}, each $c_m$ is finite,
and $c_m=0$ for all sufficiently negative $m$.

\begin{pro}\label{prop:morse-inequalities}
\textnormal{(\cite[Section~3.5]{GGM18})}
For every $m\in\Z$,
\begin{equation}\label{Morse}
  \sum_{j\leq m}(-1)^{m-j}c_j
  \geq\sum_{j\leq m}(-1)^{m-j}b_j.
\end{equation}
\end{pro}

If $\alpha$ is dynamically convex, then $c_j=b_j=0$ for
$j<n-1$ by \eqref{supp} and \eqref{b_k}.
Thus both sums in \eqref{Morse} may be taken over $n-1\leq j\leq m$.

\section{Index iteration}\label{sec:index}

\subsection{Index iteration formula}

The $\diamond$-sum of two block matrices is defined by
\[ \left(\begin{array}{lll}A_1 & B_1\\ C_1 & D_1 \end{array}\right)_{2i\times 2i}\diamond
      \left(\begin{array}{lll}A_2 & B_2\\ C_2 & D_2 \end{array}\right)_{2j\times 2j}
=\left(\begin{array}{llll}A_1 & 0 & B_1 & 0 \\
                                   0 & A_2 & 0& B_2\\
                                   C_1 & 0 & D_1 & 0 \\
                                   0 & C_2 & 0 & D_2\end{array}\right). \]
We use the basic normal forms from \cite{L00}:
\[
  N_1(\lambda,b)=
  \begin{pmatrix}\lambda&b\\0&\lambda\end{pmatrix},
  \qquad
  R(\theta)=
  \begin{pmatrix}\cos\theta&-\sin\theta\\
                  \sin\theta&\cos\theta\end{pmatrix},
\]
\[
  N_2(e^{\sqrt{-1}\theta},B)=
  \begin{pmatrix}R(\theta)&B\\0&R(\theta)\end{pmatrix},
\]
where $\lambda=\pm1$, $b\in\mathbb R$, and
$\theta\in(0,\pi)\cup(\pi,2\pi)$.
Here $B=(b_{ij})$ satisfies
$N_2(e^{\sqrt{-1}\theta},B)\in\Sp(4)$ and $b_{12}\ne b_{21}$.
This block is called \textit{non-trivial} if
$(b_{12}-b_{21})\sin\theta<0$, and \textit{trivial} if
$(b_{12}-b_{21})\sin\theta>0$.

For $M\in\Sp(2d)$, let $\sigma(M)$ denote its spectrum and set
$\nu_\omega(M)=\dim_{\mathbb{C}}\ker_{\mathbb{C}}(M-\omega I)$
for $\omega\in S^1$.
Define
\[
  \Omega(M)=\{N\in\Sp(2d)\mid
  \nu_\omega(N)=\nu_\omega(M)\text{ for every }\omega\in S^1\}.
\]
Let $\Omega^0(M)$ be the path-connected component of
$\Omega(M)$ containing $M$.

Set $\varphi(a)=\lceil a\rceil-\lfloor a\rfloor$
for $a\in\mathbb{R}$.

\begin{thm}
\label{thm:index-iteration}
\textnormal{(\cite{L00}; see also
\cite[Theorems~1.8.10 and~8.3.1]{L02})}
For every $P\in\Sp(2d)$, there exists a continuous path
$f:[0,1]\to\Omega^0(P)$ with $f(0)=P$ and
\begin{align}
  f(1)
 &=  N_1(1,1)^{\diamond p_-}\,\diamond\,I_{2p_0}\,\diamond\,N_1(1,-1)^{\diamond p_+}\nonumber \\
  &\quad \diamond\,N_1(-1,1)^{\diamond q_-}\,\diamond\,(-I_{2q_0})\,\diamond\,N_1(-1,-1)^{\diamond
q_+} \nonumber \\
  &\quad
\diamond\,N_2(e^{\sqrt{-1}\alpha_{1}},A_{1})\,\diamond\,\cdots\,\diamond\,N_2(e^{\sqrt{-1}\alpha_{r_{\ast}}},A_{r_{\ast}})
  \diamond\,N_2(e^{\sqrt{-1}\beta_{1}},B_{1})\,\diamond\,\cdots\,\nonumber \\
  &\quad
\diamond\,N_2(e^{\sqrt{-1}\beta_{r_{0}}},B_{r_{0}})\diamond\,R(\theta_1)\,\diamond\,\cdots\,\diamond\,R(\theta_r)\diamond\,H_h,\label{basic normal}
\end{align}
where the block multiplicities are nonnegative integers,
$\alpha_j,\beta_j,\theta_j\in(0,\pi)\cup(\pi,2\pi)$,
the $N_2$ blocks with angles $\alpha_j$ are non-trivial,
and those with angles $\beta_j$ are trivial.
The block $H_h\in\Sp(2h)$ is hyperbolic.
The multiplicities satisfy
\begin{equation}
p_- + p_0 + p_+ + q_- + q_0 + q_+ + r + 2r_{\ast} + 2r_0 + h = d.\label{sum_n-1}
\end{equation}

For every $\Phi\in\mathcal P_T(2d)$ with $\Phi(T)=P$
and every $m\in\N$, we have
\begin{align}
  \mu_-(\Phi^m)
 &=  m(\mu_-(\Phi)+p_-+p_0-r ) + 2\sum_{j=1}^r\left\lceil\frac{m\theta_j}{2\pi}\right\rceil - r   \nonumber
\\
  &\quad   - p_- - p_0 - \frac{1+(-1)^m}{2}(q_0+q_+)
              + 2\sum_{j=1}^{r_{\ast}}\varphi\left(\frac{m\alpha_j}{2\pi}\right) - 2r_{\ast}.
\label{3.7}
\end{align}
\end{thm}

It follows from \eqref{3.7} that
\begin{equation}\label{eq:mean-index-equality}
  |\mu_-(\Phi)-\hat\mu(\Phi)|=d
  \quad\Longleftrightarrow\quad p_-+p_0=d.
\end{equation}
\begin{equation}\label{eq:mean-index-equality1}
  |\mu_+(\Phi)-\hat\mu(\Phi)|=d
  \quad\Longleftrightarrow\quad p_++p_0=d.
\end{equation}

By \cite[Theorem~2.2]{LZ02}, if
$\Phi\in\mathcal P_T(2d)$ satisfies $\mu_-(\Phi)\geq d$, then
\begin{equation}
\mu_-(\Phi^{m+1})\ge\mu_+(\Phi^m),\ \forall\ m\in\N.
\label{ind_growth}
\end{equation}

\subsection{Common index jump theorem}

For $M\in\Sp(2d)$ and $\omega\in S^1$, let
$S_M^\pm(\omega)$ denote the splitting numbers of $M$ at $\omega$.
Their definition and the following properties can be found in
\cite[Section~9.1]{L02}.

\begin{pro}
\label{prop:splitting-numbers}
(i) $S_M^\pm(\omega)=S_N^\pm(\omega)$ for every
$N\in\Omega^0(M)$;
(ii) $S^{\pm}_M(\omega)\le \nu_\omega(M)$.

\end{pro}

Set $C(M)=\sum_{0<\theta<2\pi}S^-_M(e^{\sqrt{-1}\theta})$,
and write $\{a\}=a-\lfloor a\rfloor$ for the fractional part
of $a\in\mathbb R$.

\begin{thm}
\label{thm:common-index-jump}
\textnormal{(\cite[Theorem~4.3]{LZ02};
\cite[Theorem~3.5 and Remark~3.6]{DLW16})}
Let $\Phi_i\in\mathcal P_{T_i}(2d)$, $i=1,\ldots,q$,
be symplectic paths with $\hat\mu(\Phi_i)>0$, and set
$M_i=\Phi_i(T_i)$.
Fix a sufficiently small $\delta\in(0,1/2)$.

There exist infinitely many tuples
$(N,m_1,\ldots,m_q)\in\N^{q+1}$ such that, for all $i$,
\begin{align}
  \nu(\Phi_i) &=  \nu(\Phi_i^{2m_i\pm 1}),\label{v2m_i pm 1}
\\
  \mu_-(\Phi_i^{2m_i+1})  &=  2N+\mu_-(\Phi_i),                         \label{2m_i+1}
\\
  \mu_+(\Phi_i^{2m_i-1})  &=   2N-\mu_-(\Phi_i)-2S^+_{M_i}(1)+\nu(\Phi_i),  \label{2m_i-1}
\\
  \mu_-(\Phi_i^{2m_i})  &=  2N-(S^+_{M_i}(1)+C(M_i)-2\Delta_i), \label{2m_i}
\end{align}
where
\[
  \Delta_i=
  \sum_{\substack{0<\theta<2\pi\\
    0<\{m_i\theta/\pi\}<\delta}}
  S^-_{M_i}(e^{\sqrt{-1}\theta}).
\]
Moreover,
\begin{equation}\label{eq:central-index-upper-bound}
  \mu_+(\Phi_i^{2m_i})\leq2N+d,
  \qquad i=1,\ldots,q.
\end{equation}
For any prescribed $M_0\in\N$, the tuples can be chosen
so that $M_0\mid N$.

\end{thm}

\begin{rem}
\label{rem:index-jump-parameters}
\textup{(i)} Fix $M\in\N$ such that
$M\theta/\pi\in\Z$ whenever $\theta/\pi\in\mathbb Q$
and $e^{\sqrt{-1}\theta}\in\sigma(M_i)$ for some $i=1,\ldots,q$.
The integers $m_i$ can be chosen in the form
\begin{equation}
  m_i=\left(\left\lfloor\frac{N}{M\hat{\mu}(\Phi_i)}\right\rfloor+\chi_i\right)M,\quad\forall\  1\le i\le
q,\label{m_i}
\end{equation}
with $\chi_i\in\{0,1\}$.

\textup{(ii)} For every $e^{\sqrt{-1}\theta}\in\sigma(M_i)$,
$\{m_i\theta/\pi\}\in[0,\delta)\cup(1-\delta,1)$;
see \cite[(4.43)--(4.44)]{LZ02}.

\textup{(iii)} For any $\epsilon>0$ and $M_0\in\N$,
infinitely many such tuples satisfy $M_0\mid N$ and
\begin{equation}\label{Nv-epsilon}
  \left|
    \left\{\frac{N}{M\hat\mu(\Phi_i)}\right\}-\chi_i
  \right|<\epsilon,
  \qquad i=1,\ldots,q.
\end{equation}

\textup{(iv)} Suppose every iterate of each $\Phi_i$ is non-degenerate.
The tuples can be chosen in pairs $(N,m_1,\ldots,m_q)$ and
$(N',m_1',\ldots,m_q')$ satisfying \textup{(i)}--\textup{(iii)} and
\begin{equation}\label{eq:delta-symmetry}
  \Delta_i+\Delta_i'=C(M_i),\qquad i=1,\ldots,q.
\end{equation}
Here $\Delta_i'$ is defined by replacing $m_i$ with $m_i'$;
see \cite[Theorem~2.8]{HW22}.
\end{rem}

\section{Multiplicity and ellipticity of closed Reeb orbits}\label{sec:proofs}

Let $W\subset T^*S^n$, $n\geq3$, be a compact fiberwise
star-shaped domain containing the zero section, with contact boundary
$(M,\alpha)$.
We take coefficients in $\mathbb{Q}$.

\subsection{Symplectically degenerate maxima}

An isolated closed Reeb orbit $\gamma$ is called a \textit{symplectically degenerate maximum} (SDM) if
$\hat{\mu}(\gamma)\in2\Z$ and $\SH^{S^1}_{\hat{\mu}(\gamma)+n-1}(\gamma)\neq0$. By Section~\ref{sec:index-basic}, this definition is
independent of the symplectic trivialization of $\gamma^*\xi$. In fact, by (\ref{mu-hat{mu}}), (\ref{supp}) and (\ref{eq:mean-index-equality1}), the non-vanishing of the homology implies that $\gamma$ is totally degenerate, i.e., its linearized Poincar\'{e} map has only the eigenvalue $1$. In particular, $\hat\mu(\gamma)\in2\mathbb{Z}$ automatically holds.
We call an iteration $\gamma^k$ \textit{admissible} if
$\nu(\gamma^k)=\nu(\gamma)$.

\begin{pro}
\label{prop:sdm}
Let $\gamma\in\mathcal P(\alpha)$ be a simple isolated closed
Reeb orbit. If
\[
  \SH^{S^1}_{\hat\mu(\gamma^{k_i})+n-1}(\gamma^{k_i})\neq0
\]
for a sequence of admissible iterations $k_i\to\infty$,
then $\gamma$ is an SDM.
Moreover, the existence of a simple SDM implies that the Reeb flow
of $\alpha$ has infinitely many simple closed orbits.
\end{pro}

\begin{proof}
These facts are well-known in terms of local contact homology, see \cite[Theorem~2 and Proposition~3]{GHHM13} and \cite[Corollary~3.4]{GGM15}. It also holds for equivariant symplectic homology, as in \cite[page~570]{GG20}. For completeness, some explanation is made here.

By Theorem \ref{thm:index-iteration} and (\ref{eq:mean-index-equality1}), $\gamma$ is totally degenerate since $k_i$ is admissible. Change the trivialization along $\gamma$ so that its mean index
is zero, and use the iterated trivializations along $\gamma^{k_i}$. Then by (\ref{eq:equivariant-sh-comparison}), $\HF_{n-1}(H_{\gamma^{k_i}},p^{k_i})\neq0$, where $p$ denotes the constant orbit in Section \ref{sec:equivariant-homology}. It follows from \cite[Proposition~5.1]{GG10} that $p$ is an SDM for the Hamiltonian diffeomorphism $\varphi_{\gamma}$.

Since $\gamma$ is simple, the local contact homology of $\gamma$ is well defined and isomorphic to $\SH^{S^1}(\gamma)$ and $\HF(H_{\gamma},p)$, see \cite[Theorem~1.2]{F20} and \cite[Proposition~5.1]{HM15}.
Therefore $\gamma$ is an SDM as we defined. It is also an SDM for local contact homology as in \cite{GHHM13} and the existence of infinitely many closed orbits follows.

\end{proof}

\subsection{Proof of Theorem~\ref{multiplicity}}\label{sec:proof-multiplicity}
\leavevmode We may assume that there are finitely many simple closed Reeb
orbits, denoted by $\gamma_1,\ldots,\gamma_r$.
For each $i$, let $\Phi_i\in\mathcal{P}_{\mathcal{A}(\gamma_i)}(2n-2)$
be the symplectic path representing the linearized Reeb flow on
$\xi$ along $\gamma_i$ and set $M_i=\Phi_i(\mathcal A(\gamma_i))$.
We apply the results of Section~\ref{sec:index} with $d=n-1$.

\begin{claim}
\label{claim:1}
The mean index of $\gamma_i$ is positive.
\end{claim}

Suppose that $\hat\mu(\gamma_i)\le0$.
Dynamical convexity gives $\mu_-(\gamma_i)-\hat\mu(\gamma_i)\ge n-1$.
By \eqref{mu-hat{mu}}, equality holds, and
\eqref{eq:mean-index-equality} implies $p_{-,i}+p_{0,i}=n-1$.
Dividing \eqref{3.7} by $m$ and letting $m\to\infty$, we obtain
$\hat\mu(\gamma_i)=\mu_-(\gamma_i)+p_{-,i}+p_{0,i}
=\mu_-(\gamma_i)+n-1\ge2(n-1)>0$, a contradiction.

By Claim~\ref{claim:1}, we can apply
Theorem~\ref{thm:common-index-jump} to $\Phi_1,\ldots,\Phi_r$.
Dynamical convexity then gives
\begin{align}
  \mu_-(\gamma_i^{2m_i+1})  &\ge  2N+n-1,                         \label{DC2m_i+1}
\\
  \mu_+(\gamma_i^{2m_i-1})  &\le   2N-n+1-2S^+_{M_i}(1)+\nu(\gamma_i)\nonumber \\
  &= 2N-n+1-p_{-,i}+p_{+,i}\le2N,\label{DC2m_i-1}
\end{align}
Here the equality follows from
$S^+_{M_i}(1)=p_{-,i}+p_{0,i}$ and
$\nu(\gamma_i)=p_{-,i}+2p_{0,i}+p_{+,i}$.

Next we divide the proof into three steps.

\begin{step}
\label{step:1}
The existence of $\lfloor\frac{n}{2}\rfloor-1$ closed Reeb orbits $\gamma_1,\cdots,\gamma_{\lfloor\frac{n}{2}\rfloor-1}$.
\end{step}

This lower bound was already proved in
\cite[Theorem~1.6]{GG20}; see also \cite{W13} for closed geodesics
on Finsler spheres. We recall the argument for completeness.

Choose $N=j(n-1)$, with $j\in\mathbb N$, as allowed
by Theorem~\ref{thm:common-index-jump}.
Proposition~\ref{prop:cotangent-homology} gives closed Reeb orbits
$\eta_0,\ldots,\eta_{n-1}$ whose actions satisfy \eqref{action<}.
Their local homology satisfies
\begin{equation}\label{y_k}
  \SH^{S^1}_{2N+2k-(n-1)}(\eta_k)\neq0,
  \qquad k=0,\ldots,n-1.
\end{equation}
 For $ k\le n-2$, it follows from \eqref{ind_growth} and \eqref{DC2m_i+1} that
$\mu_-(\gamma_i^{2m_i+m})\ge \mu_-(\gamma_i^{2m_i+1})> 2N+2k-(n-1)$ for all $m\in \N$, so
\begin{equation*}
\SH^{S^1}_{2N+2k-(n-1)}(\gamma_i^{2m_i+ m})=0
\end{equation*}
  by \eqref{supp}. Similarly, for
$k\ge\lfloor\frac{n+1}{2}\rfloor$, \eqref{ind_growth}
and \eqref{DC2m_i-1} give
$\mu_+(\gamma_i^{2m_i-m})\le\mu_+(\gamma_i^{2m_i-1})<2N+2k-(n-1)$
whenever $1\le m\le2m_i-1$. Hence
\begin{equation*}
\SH^{S^1}_{2N+2k-(n-1)}(\gamma_i^{2m_i- m})=0
\end{equation*}
 by \eqref{supp}. Thus, for
$\lfloor\frac{n+1}{2}\rfloor\le k\le n-2$, the only iterate
of $\gamma_i$ that can contribute to degree $2N+2k-(n-1)$ is
$\gamma_i^{2m_i}$. By \eqref{y_k},
$\eta_k=\gamma_{i_k}^{2m_{i_k}}$ for some $i_k$.

By the strict action inequalities \eqref{action<}, these
iterates are distinct, and hence so are their underlying simple orbits.
Thus we have obtained $\lfloor\frac n2\rfloor-1$ simple closed
Reeb orbits.

\begin{step}
\label{step:2}
The existence of the $\lfloor\frac{n}{2}\rfloor$th closed Reeb orbit.
\end{step}

\begin{claim}
\label{claim:2}
There are infinitely many tuples $(N,m_1,\ldots,m_r)$
in Theorem~\ref{thm:common-index-jump}, with $N$ divisible by $n-1$,
such that
\[
  \sum_{i=1}^r2m_i\hat\chi(\gamma_i)=2N\chi_+.
\]
\end{claim}

\begin{proof}
We follow the argument used for Finsler spheres in \cite{W08}. Choose $M$ as in Remark~\ref{rem:index-jump-parameters}(i)
so that each $L_i$ in (\ref{hat chi finite sum}) and $n-1$ divide $M$, then \eqref{m_i} gives $L_i\mid m_i$ for every $i$. So
$\sum_{i=1}^r2m_i\hat\chi(\gamma_i)$ and $2N\chi_+$ are both integers by (\ref{hat chi finite sum}) and (\ref{chi_+}).

Choose $\epsilon>0$ so small that
$2M\epsilon\sum_{i=1}^r|\hat\chi(\gamma_i)|<1$,
and apply \eqref{Nv-epsilon}.
Using \eqref{m_i} and the resonance identity \eqref{MII}, we obtain
\[
\begin{aligned}
  \left|\sum_{i=1}^r2m_i\hat\chi(\gamma_i)-2N\chi_+\right|
  &=\left|2M\sum_{i=1}^r
     \left(\chi_i-\left\{\frac{N}{M\hat\mu(\gamma_i)}\right\}\right)
     \hat\chi(\gamma_i)\right|\\
  &\leq2M\epsilon\sum_{i=1}^r|\hat\chi(\gamma_i)|<1.
\end{aligned}
\]
The difference is an integer, so it is zero.
\end{proof}

Choose a tuple as in Claim~\ref{claim:2}.
We first prove that there exists a $\gamma_{i_0}$ satisfying
$\SH^{S^1}_{2N+n-1}(\gamma_{i_0}^{2m_{i_0}})\neq 0$.
Suppose instead that
\begin{equation}\label{eq:top-degree-vanishing}
  \SH^{S^1}_{2N+n-1}(\gamma_i^{2m_i})=0
  \qquad\text{for every }i=1,\ldots,r.
\end{equation}
For $k<2m_i$, \eqref{ind_growth} and \eqref{DC2m_i-1}
give $\mu_+(\gamma_i^k)\leq2N$.
For $k=2m_i$, \eqref{eq:central-index-upper-bound} and
\eqref{eq:top-degree-vanishing} exclude contributions in degrees
at least $2N+n-1$.
For $k>2m_i$, \eqref{ind_growth} and \eqref{DC2m_i+1}
give $\mu_-(\gamma_i^k)\geq2N+n-1$.
Since dynamical convexity excludes contributions below degree $n-1$,
the support inclusion \eqref{supp} gives
\begin{equation}\label{eq:truncated-euler-sum}
  \sum_{j=n-1}^{2N+n-2}(-1)^j c_j
  =\sum_{i=1}^r\sum_{k=1}^{2m_i}\chi(\gamma_i^k).
\end{equation}
By Claim~\ref{claim:2} and the Morse inequality \eqref{Morse},
we obtain a contradiction:
\begin{align*}
  (-1)^n2N\chi_+
  &=(-1)^n\sum_{i=1}^r2m_i\hat\chi(\gamma_i)\\
  &=(-1)^n\sum_{i=1}^r\sum_{k=1}^{2m_i}\chi(\gamma_i^k)\\
  &=c_{2N+n-2}-c_{2N+n-3}+\cdots\pm c_{n-1}\\
  &\geq b_{2N+n-2}-b_{2N+n-3}+\cdots\pm b_{n-1}\\
  &=(-1)^n2N\chi_++1,
\end{align*}
where the second equality follows from \eqref{hat chi finite sum},
the periodicity of $\chi(\gamma_i^k)$, and $L_i\mid m_i$;
the third follows from \eqref{eq:truncated-euler-sum},
and the last from \eqref{b_k}.

By \eqref{supp} and \eqref{eq:central-index-upper-bound},
$\mu_+(\gamma_{i_0}^{2m_{i_0}})=2N+n-1$.
Then Theorem~\ref{endpoint support} gives
\begin{equation}
\SH^{S^1}_*(\gamma_{i_0}^{2m_{i_0}})=0\quad \textnormal{for}\quad *\neq
2N+n-1,\label{one side nonzero}
\end{equation}
which implies, by \eqref{y_k}, that $\gamma_{i_0}$ is distinct
from the simple orbits found in Step~\ref{step:1}.
After relabeling, set $i_0=\lfloor\frac n2\rfloor$.
This completes Step~\ref{step:2}.

\Needspace{8\baselineskip}
\begin{step}
\label{step:3}
The existence of one more closed Reeb orbit when $n$ is odd.
\end{step}

Consider the infinitely many tuples in
Claim~\ref{claim:2}.
Suppose that, for one of these tuples, the orbit
$\eta_{\lfloor\frac{n+1}{2}\rfloor-1}$ in \eqref{y_k} equals
$\gamma_i^{2m_i}$ for some $i$.
The argument of Step~\ref{step:1} then gives
$\lfloor\frac n2\rfloor$ distinct simple orbits, and
Step~\ref{step:2} gives one additional orbit.
This proves the required bound when $n$ is odd.

Otherwise, for every such tuple,
$\eta_{\lfloor\frac{n+1}{2}\rfloor-1}=\gamma_i^{2m_i-m_0}$
for some $i$ and some $1\leq m_0\leq2m_i-1$.
Both $i$ and $m_0$ may depend on the tuple.
Since $n$ is odd, \eqref{y_k} gives
\begin{equation}
\SH^{S^1}_{2N}(\gamma_{i}^{2m_i-m_0})\neq 0\label{S3_aspt}
\end{equation}
By \eqref{supp} and \eqref{ind_growth},
\begin{equation}
2N\le\mu_+(\gamma_{i}^{2m_i-m_0})\le\mu_+(\gamma_{i}^{2m_i-1}).\label{mu_+(2mi-m0)}
\end{equation}
 Due to \eqref{DC2m_i-1}, we have the equalities
\begin{equation}
\mu_+(\gamma_{i}^{2m_i-1})=2N\label{mu_+=2N}
\end{equation}
 and
\begin{equation}
p_{+,i}=n-1.\label{p_+}
\end{equation}
 Then
\begin{equation*}
\mu_+(\gamma_{i}^{2m_i-m})\le\mu_-(\gamma_{i}^{2m_i-1})<\mu_+(\gamma_{i}^{2m_i-1})=2N
\end{equation*}
for $2\le m\le2m_i-1$.
Combining this inequality with \eqref{mu_+(2mi-m0)} gives
$m_0=1$.
\begin{claim}\label{claim: 3} \it There holds $|\hat{\mu}(\gamma_i^{2m_i})-2N|<\epsilon$ for sufficiently small $\epsilon$.\end{claim}
\begin{proof} The inequality also appeared in the {\it index recurrence theorem} of \cite{GG20}. We give an explanation in terms of the common index jump theorem. It follows from (\ref{m_i}) that
\begin{eqnarray} \hat{\mu}(\gamma_i^{2m_i})&=&2m_i\hat{\mu}(\gamma_i)=2\left(\left[\frac{N}{M\hat{\mu}(\Phi_i)}\right]+\chi_i\right)M\hat{\mu}(\gamma_i)\nonumber\\
 &=&2\left(\frac{N}{M\hat{\mu}(\Phi_i)}-\left\{\frac{N}{M\hat{\mu}(\Phi_i)}\right\}+\chi_i\right)M\hat{\mu}(\gamma_i)\nonumber\\
 &=&2N+2\left(-\left\{\frac{N}{M\hat{\mu}(\Phi_i)}\right\}+\chi_i\right)M\hat{\mu}(\gamma_i)\nonumber\end{eqnarray}
Now the result follows from (\ref{Nv-epsilon}).\end{proof}
Now we can obtain that $\gamma_i$ is an SDM for some $i$ and get a contradiction. Note that (\ref{2m_i-1}), (\ref{DC2m_i-1}) and (\ref{mu_+=2N}) imply $\mu_-(\gamma_i)=n-1$. Then by (\ref{3.7}) and (\ref{p_+}), $\hat{\mu}(\gamma_i)=n-1$ and $\hat{\mu}(\gamma_i^{2m_i})=2N$ by Claim \ref{claim: 3}. Therefore, $\hat{\mu}(\gamma_i^{2m_i-1})=2N-(n-1)$. Combined with (\ref{S3_aspt}), we have \begin{equation}
     SH^{S^1}_{\hat{\mu}(\gamma_i^{2m_i-1})+n-1}(\gamma_i^{2m_i-1})\neq0. \label{ad_neq0}
\end{equation}

Since there are only finitely many simple orbits,
some fixed $\gamma_i$ satisfies \eqref{ad_neq0} for infinitely many tuples.
By \eqref{v2m_i pm 1}, the iterates $\gamma_i^{2m_i-1}$ are admissible.
By \eqref{m_i}, $2m_i-1\to\infty$ along a subsequence
with $N\to\infty$.
Proposition~\ref{prop:sdm} then implies that $\gamma_i$ is an SDM
and that the Reeb flow has infinitely many simple closed orbits,
a contradiction!

Steps~\ref{step:1}--\ref{step:3} complete the proof
of Theorem~\ref{multiplicity}.

\subsection{Proof of Theorem~\ref{irrationally elliptic}}\label{sec:proof-ellipticity}

\leavevmode We use the notation of
Section~\ref{sec:proof-multiplicity}.
Choose two tuples $(N,m_1,\ldots,m_r)$ and
$(N',m_1',\ldots,m_r')$ as in
Remark~\partref{rem:index-jump-parameters}{iv}.
Both tuples can be chosen to satisfy the divisibility conditions
and the identity in Claim~\ref{claim:2}.

Step~\ref{step:2} gives a simple closed Reeb orbit
$\gamma_{i_1}$ with
$\SH^{S^1}_{2N+n-1}(\gamma_{i_1}^{2m_{i_1}})\neq0$.
Since $\alpha$ is non-degenerate, \eqref{supp} gives
\begin{equation}
\mu_-(\gamma_{i_1}^{2m_{i_1}})= 2N +n-1.\label{4.6aa}
\end{equation}
By \eqref{2m_i} and \eqref{4.6aa}, we have
\begin{align}
  n-1 &= -(S^+_{M_{i_1}}(1)+C(M_{i_1})-2\Delta_{i_1})\nonumber \\
  &= -(\sum_{0<\{\frac{m_{i_1}\theta}{\pi}\}<\delta}S^-_{M_{i_1}}(e^{\sqrt{-1}\theta})+
\sum_{\{\frac{m_{i_1}\theta}{\pi}\}>1-\delta}S^-_{M_{i_1}}(e^{\sqrt{-1}\theta})-2\sum_{0<\{\frac{m_{i_1}\theta}{\pi}\}<\delta}S^-_{M_{i_1}}(e^{\sqrt{-1}\theta}))\nonumber
\\
  &= \sum_{0<\{\frac{m_{i_1}\theta}{\pi}\}<\delta}S^-_{M_{i_1}}(e^{\sqrt{-1}\theta})-
\sum_{\{\frac{m_{i_1}\theta}{\pi}\}>1-\delta}S^-_{M_{i_1}}(e^{\sqrt{-1}\theta})\nonumber \\
  &\le \sum_{0<\{\frac{m_{i_1}\theta}{\pi}\}<\delta}\nu_{e^{\sqrt{-1}\theta}}(M_{i_1})\nonumber \\
  &\le r_{i_1}+r_{*,i_1}+r_{0,i_1}\le n-1,\label{equ}
\end{align}
Since $\alpha$ is non-degenerate, $S^+_{M_{i_1}}(1)=0$.
The first inequality in \eqref{equ} follows from
Proposition~\partref{prop:splitting-numbers}{ii}.
By Remark~\partref{rem:index-jump-parameters}{ii},
$\sum_{0<\{\frac{m_{i_1}\theta}{\pi}\}<\delta}\nu_{e^{\sqrt{-1}\theta}}(M_{i_1})$ is taken over exactly
one eigenvalue from each conjugate pair on the unit circle.
Each rotation block and each $N_2$ block in
\eqref{basic normal} therefore contributes one,
which gives the second inequality.
The last inequality follows from \eqref{sum_n-1}.
Since all inequalities in \eqref{equ} are equalities,
\eqref{sum_n-1} implies
\begin{equation}
  r_{i_1}=n-1.
  \label{homotopy irrationally elliptic}
\end{equation}
Thus $M_{i_1}$ is connected within $\Omega^0(M_{i_1})$
to $R(\theta_1)\,\diamond\,\cdots\,\diamond\,R(\theta_{n-1})$ and \begin{equation}
    0<\left\{\frac{m_{i_1}\theta_j}{\pi}\right\}<\delta.\label{<delta}
\end{equation}

To obtain a second orbit, we use \eqref{eq:delta-symmetry}
and \eqref{2m_i} to compute
\begin{align}
  \mu_-(\gamma_{i_1}^{2m_{i_1}'}) &= 2N'-C(M_{i_1})+2\Delta_{i_1}'\nonumber \\
  &= 2N'+C(M_{i_1})-2\Delta_{i_1}\nonumber \\
  &= 2N'-(n-1)<2N'+(n-1),\label{2m_i_1'}
\end{align}
where the last equality is due to \eqref{equ}.

Repeating the argument leading to
\eqref{4.6aa}--\eqref{homotopy irrationally elliptic}
for the primed tuple gives a simple closed Reeb orbit $\gamma_{i_2}$
such that
\begin{equation}
\mu_-(\gamma_{i_2}^{2m_{i_2}'})= 2N' +n-1\label{2m_i_2}
\end{equation}
 and its linearized Poincar\'{e} map $M_{i_2}$ is connected to
$R(\varphi_1)\,\diamond\,\cdots\,\diamond\,R(\varphi_{n-1})$ within $\Omega^0(M_{i_2})$, where
$\frac{\varphi_j}{\pi}\notin \Q, 1\le j\le n-1$. It follows from \eqref{2m_i_1'} and \eqref{2m_i_2}
that $i_1\neq i_2$.
Hence, we have obtained two simple contractible closed Reeb orbits with linearized Poincar\'{e} maps connected to a direct sum of irrational rotations within the homotopy set.

Now we prove the irrational ellipticity of $\gamma_{i_1},\gamma_{i_2}$.
For each $\theta\in\sigma(M_{i_1})$, it follows from \cite[Theorem~1.6.11]{L02} that $M_{i_1}$ is symplectically similar to $S_1 \diamond \cdots \diamond S_{m_1} \diamond S_0$, where $S_0 \in \mathrm{Sp}(2k_0)\ (k_0 \ge 0)$ does not have eigenvalue $e^{\sqrt{-1}\theta}$ and $S_i \in \mathrm{Sp}(2k_i)(k_i \ge 1)$ is of the normal form $N_{k_i}(e^{\pm\sqrt{-1}\theta}, b_i)$ defined in \cite[Section~1.6]{L02}, $\nu_{e^{\pm\sqrt{-1}\theta}}(S_i)\in\{1,2\}$.

Then by Cases 3 and 4 in \cite[Section~1.8]{L02}, if $k_i \ge 3$ for some $1 \le i \le m_1$, then $S_i$ can be connected within $\Omega^0(S_i)$ to $\tilde{S}_i$ with $e(\tilde{S}_i) < e(S_i)$. This contradicts to (\ref{homotopy irrationally elliptic}). If $k_i = 2$ for some $1 \le i \le m_1$, then by (\ref{homotopy irrationally elliptic}), $S_i$ is not a basic normal form, hence $\nu_{e^{\pm\sqrt{-1}\theta}}(S_i)=2$ by \cite[Lemma~1.9.2]{L02}. By \cite[Case~4 in Section~1.8]{L02}, $S_i$ can be connected within $\Omega^0(S_i)$ to $R(\theta) \diamond R(2\pi - \theta)$. This contradicts to (\ref{<delta}). Now $M_{i_1}$ is symplectically similar to $R(\omega)^{\diamond m_1} \diamond S_0$
with $\omega = \theta$ or $2\pi - \theta$. By the inductive argument, we eventually obtain that $M_{i_1}$ is symplectically similar to $R(\theta_1)\,\diamond\,\cdots\,\diamond\,R(\theta_{n-1})$. The same argument proves the irrational ellipticity of $x_{i_2}$.

\section*{Funding}
H. Duan was partially supported by NSFC (12631007, 12271268 and 12361141812) and Natural Science Foundation of Tianjin (25JCZDJC01030).
Z. Qi was partially supported by NSFC (125B200041). D. Xie was partially supported by NSFC (12631007 and 12271268).

\end{document}